\documentclass[letterpaper]{article}
\usepackage{aaai}
\usepackage{times}
\usepackage{helvet}
\usepackage{courier}
\usepackage{graphicx}
\usepackage{amsmath}
\usepackage{array}
\usepackage{multirow}
\usepackage{mdwmath}
\usepackage{mdwtab}
\usepackage[tight,footnotesize]{subfigure}
\usepackage{amsfonts}
\usepackage{algorithm}
\usepackage{algorithmic}
\usepackage{balance}

\newtheorem{theorem}{Theorem}

\begin{document}
%
\title{Generalized Higher-Order Tensor Decomposition via Parallel ADMM}
\author{Fanhua Shang$^{1}$, Yuanyuan Liu$^{2}$, James Cheng$^{1}$\\
$^{1}$Department of Computer Science and Engineering, The Chinese University of Hong Kong\\
$^{2}$Department of Systems Engineering and Engineering Management, The Chinese University of Hong Kong\\
$\{$fhshang, jcheng$\}$@cse.cuhk.edu.hk; yyliu@se.cuhk.edu.hk}
\maketitle

\begin{abstract}
\begin{quote}
Higher-order tensors are becoming prevalent in many scientific areas such as computer vision, social network analysis, data mining and neuroscience. Traditional tensor decomposition approaches face three major challenges: model selecting, gross corruptions and computational efficiency. To address these problems, we first propose a parallel trace norm regularized tensor decomposition method, and formulate it as a convex optimization problem. This method does not require the rank of each mode to be specified beforehand, and can automatically determine the number of factors in each mode through our optimization scheme. By considering the low-rank structure of the observed tensor, we analyze the equivalent relationship of the trace norm between a low-rank tensor and its core tensor. Then, we cast a non-convex tensor decomposition model into a weighted combination of multiple much smaller-scale matrix trace norm minimization. Finally, we develop two parallel alternating direction methods of multipliers (ADMM) to solve our problems. Experimental results verify that our regularized formulation is effective, and our methods are robust to noise or outliers.
\end{quote}
\end{abstract}

\section{Introduction}
The term tensor used in the context of this paper refers to a multi-dimensional array, also known as a multi-way or multi-mode array. For example, if $\mathcal{X}\in\mathbb{R}^{{I_{1}}\times{I_{2}}\times{I_{3}}}$, then we say $\mathcal{X}$ is a third-order tensor, where order is the number of ways or modes of the tensor. Thus, vectors and matrices are first-order and second-order tensors, respectively. Higher-order tensors arise in a wide variety of application areas, such as machine learning (Tomioka and Suzuki, 2013; Signoretto et al., 2014), computer vision (Liu et al., 2009), data mining (Yilmaz et al., 2011; Morup, 2011; Narita et al., 2012; Liu et al., 2014), numerical linear algebra (Lathauwer et al., 2000a; 2000b), and so on. Especially, with the rapid development of modern computer technology in recent years, tensors are becoming increasingly ubiquitous such as multi-channel images and videos, and have become increasingly popular (Kolda and Bader, 2009). When working with high-order tensor data, various new computational challenges arise due to the exponential increase in time and memory space complexity when the number of orders increases. This is called the curse of dimensionality. In practice, the underlying tensor is often low-rank, even though the actual data may not be due to noise or arbitrary errors. Essentially, the major component contained in the given tensor is often governed by a relatively small number of latent factors.

One standard tool to alleviate the curse is tensor decomposition. Decomposition of high-order tensors into a small number of factors has been one of the main tasks in multi-way data analysis, and commonly takes two forms: Tucker decomposition (Tucker, 1966) and CANDECOMP/PARAFAC (CP) (Harshman, 1970) decomposition. There are extensive studies in the literature for finding the Tucker decomposition and the CP decomposition for higher-order tensors (Kolda and Bader, 2009). In those tensor decomposition methods, their goal is to (approximately) reconstruct the input tensor as a sum of simpler components with the hope that these simpler components would reveal the latent structure of the data. However, existing tensor decomposition methods face three major challenges: rank selection, outliers and gross corruptions, and computational efficiency. Since the Tucker and CP decomposition methods are based on least-squares approximation, they are also very sensitive to outliers and gross corruptions (Goldfarb and Qin, 2014). In addition, the performance of those methods is usually sensitive to the given ranks of the involved tensor (Liu et al., 2009). To address the problems, we propose two robust and parallel higher-order tensor decomposition methods with trace norm regularization.

Recently, much attention has been drawn to the low-rank tensor recovery problem, which arises in a number of real-word applications, such as 3D image recovery, video inpainting, hyperspectral data recovery, and face reconstruction. Compared with matrix-based analysis methods, tensor-based multi-linear data analysis has shown that tensor models are capable of taking full advantage of the high-order structure to provide better understanding and more precision. The key idea of low-rank tensor completion and recovery methods is to employ matrix trace norm minimization (also known as the nuclear norm, which is the convex surrogate of the rank of the involved matrix). In addition, there are some theoretical developments that guarantee reconstruction of a low-rank tensor from partial measurements or grossly corrupted observations via solving the trace norm minimization problem under some reasonable conditions (Tomioka et al., 2011; Shi et al., 2013). Motivated by the recent progress in tensor completion and recovery, one goal of this paper is to extend the trace norm regularization to robust higher-order tensor decomposition.

Different from existing tensor decomposition methods, we first propose a parallel trace norm regularized tensor decomposition method, which can automatically determine the number of factors in each mode through our optimization scheme. In other words, this method does not require the rank of each mode to be specified beforehand. In addition, by considering the low-rank structure of the observed tensor and further improving the scalability of our convex method, we analyze the equivalent relationship of the trace norm between a low-rank tensor and its core tensor. Then, we cast the non-convex trace norm regularized higher-order orthogonal iteration model into a weighted combination of multiple much-smaller-scale matrix trace norm minimization. Moreover, we design two parallel alternating direction methods of multipliers (ADMM) to solve the proposed problems, which are shown to be fast, insensitive to initialization and robust to noise and/or outliers with extensive experiments.

\section{Notations and Related Work}
We first introduce the notations, and more details can be seen in Kolda and Bader (2009). An $N$th-order tensor is denoted by a calligraphic letter, e.g., $\mathcal{T}\in\mathbb{R}^{{I_{1}}\times{I_{2}}\times\cdots\times{I_{N}}}$, and its entries are denoted by $t_{{i_{1}}\cdots{i_{n}}\cdots{i_{N}}}$, where $i_{n}\in\{1,\cdots,I_{n}\}$ for $1\leq n\leq N$. Fibers are the higher-order analogue of matrix rows and columns. The mode-$n$ fibers are vectors $\mathrm{t}_{{i_{1}}\cdots{i_{n-1}}{i_{n+1}}\cdots{i_{N}}}$ that are obtained by fixing the values of $\{i_{1},\cdots,i_{N}\}\backslash{i_{n}}$.

The mode-$n$ unfolding, also known as matricization, of an $N$th-order tensor $\mathcal{T}$ is denoted by $\mathcal{T}_{(n)}\in\mathbb{R}^{{{I_{n}}\times{\Pi_{j\neq{n}}}{I_{j}}}}$ and arranges the mode-$n$ fibers to be the columns of the resulting matrix $\mathcal{T}_{(n)}$ such that the mode-$n$ fiber becomes the row index and all other $(N-1)$ modes become the column indices. The tensor element $(i_{1},i_{2},\cdots,i_{N})$ is mapped to the matrix element $(i_{n}, j)$, where
\begin{displaymath}
j=1+\sum^{N}_{k=1,k\neq n}(i_{k}-1)J_{k}\;\;\textup{with}\;\;J_{k}=\prod^{k-1}_{m=1, m\neq n}I_{m}.
\end{displaymath}

The inner product of two same-sized tensors $\mathcal{A}\in \mathbb{R}^{{I_{1}}\times{I_{2}}\times}$ $^{\cdots\times{I_{N}}}$ and $\mathcal{B}\in\mathbb{R}^{{I_{1}}\times{I_{2}}\times\cdots\times{I_{N}}}$ is the sum of the product of their entries, $<{\mathcal{A},\mathcal{B}}>=\sum_{i_{1},\cdots{i_{N}}}a_{{i_{1}}\cdots{i_{N}}}b_{{i_{1}}\cdots{i_{N}}}$. The Frobenius norm of an $N$th-order $\mathcal{T}$ is defined as:
\begin{displaymath}
{\|\mathcal{T}\|}_{F}:=\sqrt{\sum^{I_{1}}_{i_{1}=1}\cdots\sum^{I_{N}}_{i_{N}=1}t^{2}_{{i_{1}}\cdots{i_{N}}}}.
\end{displaymath}

Let $A$ and $B$ be two matrices of size $m\times n$ and $p\times q$, respectively. The Kronecker product of two matrices $A$ and $B$, denoted by $A\otimes B$, is an $mp\times{nq}$ matrix given by:
\begin{displaymath}
A\otimes B=[a_{ij}B]_{mp\times{nq}}.
\end{displaymath}
The $n$-mode matrix product of a tensor $\mathcal{T}$ with a matrix $U\in \mathbb{R}^{\emph{J}\times{I_{n}}}$, denoted by $\mathcal{T}\times_{n}U\in \mathbb{R}^{{I_{1}}\times\cdots{I_{n-1}}\times J\times{I_{n+1}}\times\cdots\times{I_{N}}}$, is defined as:
\begin{displaymath}
(\mathcal{T}\times_{n}U)_{{i_{1}}\cdots{i_{n-1}}j{i_{n+1}}\cdots{i_{N}}}=\sum^{I_{n}}_{i_{n}=1}t_{{i_{1}}{i_{2}}\cdots{i_{N}}}u_{j{i_{n}}}.
\end{displaymath}

\subsection{Tensor Decomposition}
We will review two popular models for tensor decomposition, i.e., the Tucker decomposition and the CANDECOMP/PARAFAC (CP) decomposition. It is well known that finding the CP decomposition with the minimum tensor rank is a hard problem, and there is no straightforward algorithm for computing the rank for higher-order tensors (Hillar and Lim, 2013). The Tucker decomposition decomposes a given tensor $\mathcal{T}$ into a core tensor multiplied by a factor matrix along each mode as follows:
\begin{equation}
\mathcal{T}=\mathcal{G}\times_{1}U_{1}\times_{2}\cdots\times_{N}U_{N},
\end{equation}
where $U_{n}\in \mathbb{R}^{{I_{n}}\times{R_{n}}}$ are the factor matrices, which can be thought of as the principal components in each mode, and the entries of the core tensor $\mathcal{G}\in\mathbb{R}^{{R_{1}}\times{R_{2}}\cdots\times{R_{N}}}$ show the level of interaction between the different components. Since the decomposition rank $R_{n}\;(n=1,\cdots, N)$ is in general much smaller than $I_{n}\;(n=1,\cdots, N)$, the core tensor $\mathcal{G}$ can be thought of as a low-rank version of $\mathcal{T}$ (e.g., the Tucker decomposition of a third-order tensor is illustrated in Figure 1). In this sense, the storage of the Tucker decomposition form can be significantly smaller than that of the original tensor. Moreover, unlike the rank of the tensor, $R_{n}$, i.e., the mode-$n$ rank $(n=1,\cdots, N)$, is the rank of the mode-$n$ unfolding, and is clearly computable. Hence, we are particularly interested in extending the Tucker decomposition for tensor analysis.
\begin{figure}
\centering
\includegraphics [width=0.75\linewidth] {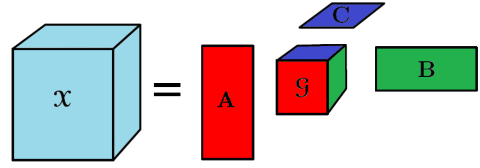}
\caption{Illustration of the Tucker decomposition of a third-order tensor.}
\label{fig_sim}
\end{figure}

If the factor matrices of the Tucker decomposition are constrained orthogonal, the decomposition form is referred to as the higher-order singular value decomposition (HOSVD, Lathauwer et al., 2000a) or higher-order orthogonal iteration (HOOI, Lathauwer et al., 2000b), where the latter leads to the estimation of best rank-$(R_{1},R_{2},\ldots,R_{N})$ approximations while the truncation of HOSVD may achieve a good rank-$(R_{1},R_{2},\ldots,R_{N})$ approximation but in general not the best possible one (Lathauwer et al., 2000b). Therefore, HOOI is the most widely used tensor decomposition method (Kolda and Bader, 2009), and takes the following form
\begin{equation}
\begin{split}
&\min_{\mathcal{G}, U_{n}}\,\|{\mathcal{T}-\mathcal{G}\times_{1}U_{1}\times_{2}\cdots\times_{N}U_{N}}\|^{2}_{F},\\
&\;\textup{s.t.},\,U^{T}_{n}U_{n}=I_{R_{n}}, n=1,\cdots, N.
\end{split}
\end{equation}

\section{Convex Tensor Decomposition Method}
\subsection{Convex Tensor Decomposition Model}
Given a tensor $\mathcal{T}$, our goal is to find a low-rank tensor $\mathcal{X}$, in order to minimize the Frobenius norm of their difference as follows:
\begin{equation}
\min_{\mathcal{X}}\,\frac{1}{2}{\|\mathcal{X}-\mathcal{T}\|^{2}_{F}}.
\end{equation}

Different from the matrix case, the low-rank tensor estimation problem (3) is in general hard to solve (Narita et al., 2012). Following the progress in tensor completion, we cast it into a (weighted) trace norm minimization problem:
\begin{equation}
\min_{\mathcal{X}}\,\sum^{N}_{n=1}\|\mathcal{X}_{(n)}\|_{\textup{tr}}+\frac{\lambda}{2}\|\mathcal{X}-\mathcal{T}\|^{2}_{F},
\end{equation}
where $\|\mathcal{X}_{(n)}\|_{\textup{tr}}$ denotes the trace norm of the unfolded matrix $\mathcal{X}_{(n)}$, i.e., the sum of its singular values, $\lambda>0$ is a regularization parameter. For handling the unbalanced target tensor, we briefly introduce the preselected weights $\alpha_{n}\geq 0$ (satisfying $\sum_{n}\alpha_{n}=1$) to the trace norm term in (4).

\subsection{Parallel Optimization Algorithm}
Due to the interdependent matrix trace norm terms, the proposed tensor decomposition model (4) is very difficult to solve. Thus, we introduce some auxiliary variables $\mathcal{M}_{n}$ into the model (4) and reformulate it into the following equivalent form:
\begin{equation}
\begin{split}
&\min_{\mathcal{X},\{\mathcal{M}_{n}\}}\,\sum^{N}_{n=1}\|\mathcal{M}_{n,(n)}\|_{\textup{tr}}+\frac{\lambda}{2}\|\mathcal{X}-\mathcal{T}\|^{2}_{F},\\
&\;\;\textup{s.t.},\,\mathcal{M}_{n}=\mathcal{X},\,n=1,\cdots, N.
\end{split}
\end{equation}

In the following, we will design a parallel alternating direction method of multipliers (ADMM) for solving the problem (5). The algorithm decomposes the original problem into smaller subproblems and solves them individually in parallel at each iteration. The parallel ADMM for the problem (5) is derived by minimizing the augmented Lagrange function $\mathcal{L}_{\mu}$ with respect to $(\mathcal{X},\;\{\mathcal{M}_{n}\})$ in parallel, and then updating the multiplier tensor $\mathcal{Y}_{n}$. But the parallelization is likely to diverge, therefore, it is necessary to modify the common algorithm in a certain way to guarantee its convergence. Several variants of parallel ADMM have been proposed in He (2009) and Deng et al., (2013) by taking additional correction steps at every iteration. Similar to Deng et al., (2013), we add some proximal terms to each subproblems and attach a relaxation parameter $\gamma>0$ for the update of the penalty parameter $\mu>0$.

\subsubsection{Computing $\{\mathcal{M}^{k+1}_{n},\;n=1,\;\cdots,\;N\}$:}
By keeping all other variables fixed, the optimal $\mathcal{M}^{k}_{n}$ is the solution to the following problem:
\begin{equation}
\begin{split}
\min_{\mathcal{M}_{n}}\,\|\mathcal{M}_{n,\;(n)}\|_{\textup{tr}}&+\frac{\mu}{2}\|\mathcal{M}_{n}-\mathcal{X}^{k}+\mathcal{Y}^{k}_{n}/\mu\|^{2}_{F}\\
&+\frac{\tau_{n}}{2}\|\mathcal{M}_{n}-\mathcal{M}^{k}_{n}\|^{2}_{F},
\end{split}
\end{equation}
where $\tau_{n}>0$ is a constant for the proximal term $\|\mathcal{M}_{n}-\mathcal{M}^{k}_{n}\|^{2}_{F}$. Following (Cai et al., 2010), we obtain a closed-form solution to the problem (6) via the following proximal operator of the trace norm:
\begin{equation}
\mathcal{M}^{k+1}_{n}= \textup{refold}(\textup{prox}_{\frac{1}{\mu+\tau_{n}}}(\frac{\mu\mathcal{X}^{k}_{(n)}-\mathcal{Y}^{k}_{(n)}+\tau_{n}\mathcal{M}^{k}_{(n)}}{\mu+\tau_{n}})),
\end{equation}
where $\textup{refold}(\cdot)$ denotes the refolding of the matrix into a tensor, i.e., the reverse process of unfolding.

\subsubsection{Computing $\mathcal{X}^{k+1}$:}
The optimal $\mathcal{X}^{k+1}$ with all other variables fixed is the solution to the following problem:
\begin{equation}
\begin{split}
\min_{\mathcal{X}}\,\frac{\lambda}{2}\|\mathcal{X}-\mathcal{T}\|^{2}_{F}&+\sum^{N}_{n=1}\frac{\mu}{2}\|\mathcal{M}^{k}_{n}-\mathcal{X}+\mathcal{Y}^{k}_{n}/\mu\|^{2}_{F}\\
&+\frac{\tau_{N+1}}{2}\|\mathcal{X}-\mathcal{X}^{k}\|^{2}_{F},
\end{split}
\end{equation}
where $\tau_{N+1}>0$ is a constant for the proximal term $\|\mathcal{X}-\mathcal{X}^{k}\|^{2}_{F}$. Since the problem (8) is a smooth convex optimization problem, it is easy to show that the optimal solution to (8) is given by
\begin{equation}
\mathcal{X}^{k+1}=\frac{\sum^{N}_{n=1}(\mu\mathcal{M}^{k}_{n}+\mathcal{Y}^{k}_{n})+\lambda\mathcal{T}+\tau_{N+1}\mathcal{X}^{k}}{N\mu+\lambda+\tau_{N+1}}.
\end{equation}

Based on the description above, we develop a parallel ADMM algorithm for the convex tensor decomposition (CTD) problem (5), as outlined in \textbf{Algorithm 1}. In Algorithm 1, we use a Jacobi-type scheme to update $(N+1)$ blocks of variables, $\{\mathcal{M}_{1}, \cdots, \mathcal{M}_{N}\}$ and $\mathcal{X}$. By the definition $f(\mathcal{M}_{1},\cdots,\mathcal{M}_{N}):=\sum^{N}_{n=1}\|\mathcal{M}_{n,(n)}\|_{\textup{tr}}$ and $g(\mathcal{X}):=\frac{\lambda}{2}{\|\mathcal{X}-\mathcal{T}\|^{2}_{F}}$, it is easy to verify that the problem (5) and Algorithm 1 satisfy the convergence conditions in Deng et al., (2013).
\begin{algorithm}[t]
\caption{Solving problem (5) via parallel ADMM}
\label{alg:Framwork2}
\renewcommand{\algorithmicrequire}{\textbf{Input:}}
\renewcommand{\algorithmicensure}{\textbf{Initialize:}}
\renewcommand{\algorithmicoutput}{\textbf{Output:}}
\begin{algorithmic}[1]
\REQUIRE Given $\mathcal{T}$, $\lambda$, and $\mu$.\\
\ENSURE $\mathcal{M}^{0}_{n}=\mathcal{X}^{0}=\mathcal{Y}^{0}_{n}=0$, $\forall n\in\{1,\cdots,N\}$.\\
\WHILE {not converged}
\FOR {$n=1, \cdots, N$}
\STATE {Update $\mathcal{M}^{k+1}_{n}$ by (7),\\
where $U_{n}S_{n}V^{T}_{n}=\frac{\mu\mathcal{X}^{k}_{(n)}-\mathcal{Y}^{k}_{(n)}+\tau_{n}\mathcal{M}^{k}_{(n)}}{\mu+\tau_{n}}$.}
\ENDFOR
\STATE{Update $\mathcal{X}^{k+1}$ by (9).}
\FOR {$n=1, \cdots, N$}
\STATE {$\mathcal{Y}^{k+1}_{n}=\mathcal{Y}^{k}_{n}+\gamma\mu(\mathcal{M}^{k+1}_{n}-\mathcal{X}^{k+1})$}.
\ENDFOR
\STATE {Check the convergence condition,\\
$\|\mathcal{X}^{k+1}-\mathcal{T}\|_{F}<\textup{Tol}$.}
\ENDWHILE
\STATE {$\mathcal{G}=\mathcal{X}^{k+1}\times_{1}U^{T}_{1}\cdots\times_{N}U^{T}_{N}$}.
\OUTPUT $\mathcal{X}^{k}$, $\mathcal{G}$, and $U_{n},\;n=1,\cdots, N$.
\end{algorithmic}
\end{algorithm}

\begin{theorem}
Let $\tau_{i}>\mu(\frac{N}{2-\gamma}-1),\,i=1,\ldots,N+1$. Then the sequence $\{\mathcal{M}^{k}_{1}, \cdots, \mathcal{M}^{k}_{N}, \mathcal{X}^{k}\}$ generated by Algorithm 1 converges to an optimal solution $\{\mathcal{M}^{\ast}_{1}, \cdots, \mathcal{M}^{\ast}_{N},\mathcal{X}^{\ast}\}$ of the problem (5). Hence, the sequence $\{\mathcal{X}^{k}\}$ converges to an optimal solution to the low-rank tensor decomposition problem (4).
\end{theorem}

Our convex Tucker decomposition method employs matrix trace norm regularization and uses the singular value decomposition (SVD) in Algorithm 1, which becomes a little slow or even not applicable for large-scale problems. Motivated by this, we will propose a scalable non-convex low-rank tensor decomposition method.

\section{Non-Convex Tensor Decomposition Method}
Several researchers (Keshavan et al., 2010; Wen et al., 2012) have provided some matrix rank estimation strategies to compute some good values $(r_{1}, r_{2},\ldots, r_{N})$ for the \emph{N}-rank of the involved tensor. Thus, we only set some relatively large integers $(R_{1}, R_{2},\ldots, R_{N})$ such that $R_{n}\geq r_{n}$, $n=1, \cdots, N$.
\begin{theorem}
Let $\mathcal{X}\in \mathbb{R}^{{I_{1}}\times{I_{2}}\times\cdots\times{I_{N}}}$ with $N$-rank$=(r_{1},r_{2}$, $\cdots,r_{N})$ and $\mathcal{G}\in\mathbb{R}^{{R_{1}}\times{R_{2}}\times\cdots\times{R_{N}}}$ satisfy $\mathcal{X}=\mathcal{G}\times_{1}U_{1}\cdots\times_{N}U_{N}$, and $U^{T}_{n}U_{n}=I_{R_{n}}$, $n=1,2,\cdots,N$, then
\begin{equation}
\|\mathcal{X}_{(n)}\|_{\textup{tr}}=\|\mathcal{G}_{(n)}\|_{\textup{tr}},\;n=1, 2, \cdots, N.
\end{equation}
\end{theorem}
The proof of Theorem 2 can be found in the supplemental material. According to the theorem, it is cleat that although the core tensor $\mathcal{G}$ of size $({R_{1}},{R_{2}},\cdots,{R_{N}})$ has much smaller sizes than the original tensor $\mathcal{X}$ (usually, $R_{n}\ll I_{n},n=1,2,\cdots,N $), their trace norm is identical. In other words, each unfolding $\mathcal{G}_{(n)}\in\mathbb{R}^{R_{n}\times\Pi_{j\neq{n}}{R_{j}}}$ of the core tensor $\mathcal{G}$ has much smaller sizes than that of the original tensor, $\mathcal{X}_{(n)}\in\mathbb{R}^{I_{n}\times\Pi_{j\neq{n}}{I_{j}}}$. Therefore, we use the trace norm of each unfolding of the core tensor to replace that of the original tensor.

\subsection{Generalized HOOI Model with Trace Norm Penalties}
According to the analysis above, our trace norm regularized HOOI model is formulated into the following form:
\begin{equation}
\begin{split}
&\min_{\mathcal{G},\{U_{n}\}}\,\sum^{N}_{n=1}{\|\mathcal{G}_{(n)}\|_{\textup{tr}}}+\frac{\lambda}{2}\| \mathcal{T}-\mathcal{G}\times_{1}U_{1}\cdots\times_{N}U_{N}\|^{2}_{F},\\
&\;\textup{s.t.},\,U^{T}_{n}U_{n}=I_{R_{n}}, n=1,2,\cdots,N.
\end{split}
\end{equation}
The core tensor trace norm regularized HOOI model (11), also called Generalized HOOI, can alleviate the SVD computational burden of large unfolded matrices involved in the convex Tucker decomposition problem (4). Furthermore, we use the trace norm regularization term in (11) to promote the robustness of the rank selection, while the original Tucker decomposition method is usually sensitive to the given ranks $({R_{1}},{R_{2}},\cdots,{R_{N}})$ (Liu et al., 2009). Due to the interdependent matrix trace norm terms, we apply the variable-splitting technique and introduce some auxiliary variables $G_{n}\in\mathbb{R}^{{{R_{n}}\times{\Pi_{j\neq{n}}}{R_{j}}}}$ into our model (11), and then reformulate the model (11) into the following equivalent form:
\begin{equation}
\begin{split}
&\min_{\mathcal{G},\{U_{n},G_{n}\}}\,\sum^{N}_{n=1}{\|G_{n}\|_{\textup{tr}}}+\frac{\lambda}{2}\| \mathcal{T}-\mathcal{G}\times_{1}U_{1}\cdots\times_{N}U_{N} \|^{2}_{F},\\
&\qquad{\textup{s.t.},\,G_{n}=\mathcal{G}_{(n)},\,U^{T}_{n}U_{n}=I_{R_{n}}, n=1,2,\cdots,N.}
\end{split}
\end{equation}

\subsection{Parallel Optimization Procedure}
In this part, we will also propose a parallel ADMM algorithm to solve the problem (12).

\subsubsection{Updating $\{\mathcal{G}^{k+1},U^{k+1}_{1},\cdots,U^{k+1}_{N}\}$:}
The optimization problem with respect to $\mathcal{G}$ and $\{U_{1},\cdots, U_{N}\}$ is formulated as follows:
\begin{equation}
\begin{split}
\mathop{\min}_{\mathcal{G},\{U_{n}\}}\sum^{N}_{n=1}\frac{\mu}{2}\|\mathcal{G}_{(n)}-G^{k}_{n}+Y^{k}_{n}/\mu\|^{2}_{F}\\
+\frac{\lambda}{2}\|\mathcal{T}-\mathcal{G}\times_{1}U_{1}\cdots\times_{N}U_{N}\|^{2}_{F}.
\end{split}
\end{equation}
Following (Lathauwer et al., 2000b), it is sufficient to determine the matrices $\{U_{1},\cdots, U_{N}\}$ for the optimization of the problem (13). For any estimate of these matrices, the optimal solution with respect to $\mathcal{G}$ is given by the following theorem.

\begin{theorem}
For given matrices $\{U_{1},\cdots,U_{N}\}$, the optimal core tensor $\mathcal{G}$ to the optimization problem (13) is given by
\begin{equation}
\begin{split}
\mathcal{G}^{k+1}=&\frac{\lambda}{\lambda+N\mu}\mathcal{T}\times_{1}(U^{k}_{1})^{T}\cdots\times_{N}(U^{k}_{N})^{T}\\
&+\frac{\mu}{\lambda+N\mu}\sum^{N}_{n=1}\textup{refold}(G^{k}_{n}-Y^{k}_{n}/\mu).
\end{split}
\end{equation}
\end{theorem}
The proof of Theorem 3 can be found in the supplemental material. In the following, we design a generalized higher order orthogonal iteration scheme for solving $\{U_{1},\cdots,U_{N}\}$ that is an alternating least squares (ALS) approach to solve the rank-$(R_{1},\ldots,R_{N})$ problem. Analogous with Theorem 4.2 in Lathauwer et al., (2000b), we first state that the minimization problem (13) can be formulated as follows.
\begin{theorem}
Assume a real $N$th-order tensor $\mathcal{T}\in\mathbb{R}^{I_{1}\times I_{2}\times\ldots I_{N}}$. Then the minimization of the problem (13) over the matrices $U_{1},\ldots,U_{N}$ having orthonormal columns is equivalent to the maximization of the following problem
\begin{equation}
h(U_{1},\ldots,U_{N})=\langle\mathcal{T},\,\mathcal{G}\times_{1}U_{1}\cdots\times_{N}U_{N}\rangle.
\end{equation}
\end{theorem}
The proof of Theorem 4 can be found in the supplemental material. According to the theorem, a generalized higher-order orthogonal iteration scheme is proposed to solve the problem (15) that solves $U_{n},\,n=1,\ldots,N$ by fixing other variables $U_{j},j\neq n $ in parallel. Imagine that the matrices $\{U_{1},\ldots,U_{n-1}$, $U_{n+1},\ldots,U_{N}\}$ are fixed and that the optimization problem (13) is thought of as a quadratic expression in the components of the matrix $U_{n}$ that is being optimized. Considering that the matrix having orthonormal columns, we have
\begin{equation}
h=\textup{trace}((U^{k+1}_{n})^{T}\mathcal{T}_{(n)}W^{T}_{n}),
\end{equation}
where $\textup{trace}(A)$ denotes the trace of the matrix $A$, and
\begin{equation}
\begin{split}
W_{n}=\mathcal{G}_{(n)}&[(U^{k}_{1})^{T}\cdots\otimes(U^{k}_{n-1})^{T}\\
&\otimes(U^{k}_{n+1})^{T}\cdots\otimes(U^{k}_{N})^{T}].
\end{split}
\end{equation}
This is actually the well-known orthogonal procrustes problem (Nick, 1995), whose global optimal solution is given by the singular value decomposition of $\mathcal{T}_{(n)}W^{T}_{n}$, i.e.,
\begin{equation}
U^{k+1}_{n}=\widehat{U_{n}}(\widehat{V_{n}})^{T},
\end{equation}
where $\widehat{U_{n}}$ and $\widehat{V_{n}}$ are obtained by SVD of $\mathcal{T}_{(n)}W^{T}_{n}$. Repeating the procedure above in parallel for different modes leads to an alternating least squares scheme for solving the maximization of the problem (15).

\subsubsection{Updating $\{G^{k+1}_{1},\cdots,G^{k+1}_{N}\}$:}
By keeping all other variables fixed, $G^{k+1}_{n}$ is updated by solving the following problem:
\begin{equation}
\mathop{\min}_{G_{n}}\|G_{n}\|_{\textup{tr}}+\frac{\mu}{2}\|G_{n}-\mathcal{G}^{k}_{(n)}-Y^{k}_{n}/\mu\|^{2}_{F}+\frac{\tau_{n}}{2}\|G_{n}-G^{k}_{n}\|^{2}_{F}.
\end{equation}
Similar to the problem (6), it is easy to obtain a closed-form solution to the problem (19):
\begin{equation}
G^{k+1}_{n}=\textup{prox}_{1/(\mu+\tau_{n})}(\frac{\mu\mathcal{G}^{k}_{(n)}+Y^{k}_{n}+\tau_{n}G^{k}_{n}}{\mu+\tau_{n}}).
\end{equation}
From (20), it is clear that only some smaller sized matrices $(\mu\mathcal{G}^{k}_{(n)}+Y^{k}_{n}+\tau_{n}G^{k}_{n})/(\mu+\tau_{n})\in\mathbb{R}^{{{R_{n}}\times{\Pi_{j\neq{n}}}{R_{j}}}}$ need to perform SVD. Thus, our non-convex trace norm regularized method has a significantly lower computational complexity with $O(\sum_{n}{{R^{2}_{n}}\times{\Pi_{j\neq{n}}}{R_{j}}})$, while the computational complexity of our convex algorithm for the problem (4) is $O(\sum_{n}{{I^{2}_{n}}\times{\Pi_{j\neq{n}}}{I_{j}}})$ at each iteration.

Based on the analysis above, we develop a parallel ADMM algorithm for solving the low-rank non-convex tensor decomposition (NCTD) problem (11), as outlined in \textbf{Algorithm 2}. Our algorithm is essentially a Jacobi-type scheme of ADMM, and is well suited for parallel and distributed computing and hence is particularly attractive for solving certain large-scale problems. Moreover, the update strategy of Gauss-Seidel version of ADMM is easily implemented. This algorithm can be accelerated by adaptively changing $\mu$. An efficient strategy (Lin et al., 2011) is to let $\mu^{0}=\mu$ (initialized in Algorithm 1 and Algorithm 2) and increase $\mu^{k}$ iteratively by $\mu^{k+1}=\rho\mu^{k}$, where $\rho\in (1.0,1.1]$ in general and $\mu^{0}$ is a small constant. Moreover, the stability and efficiency of our NCTD method can be validated in the experimental section.

\begin{algorithm}
\caption{Solving problem (11) via Parallel ADMM}
\label{alg:Framwork3}
\renewcommand{\algorithmicrequire}{\textbf{Input:}}
\renewcommand{\algorithmicensure}{\textbf{Initialize:}}
\renewcommand{\algorithmicoutput}{\textbf{Output:}}
\begin{algorithmic}[1]
\REQUIRE $\mathcal{T}$, the tensor ranks\,$({R_{1}},{R_{2}},\cdots,{R_{N}})$, and $\lambda$.\\
\ENSURE $Y^{0}_{n}=0$, $G^{0}_{n}=0$, $U^{0}_{n}=\textup{rand}(I_{n},\,R_{n})$, $\mu^{0}=10^{-4}$, $\mu_{max}=10^{10}$, $\rho=1.05$, and $\textup{Tol}=10^{-5}$.\\
\WHILE {not converged}
\STATE {Update $\mathcal{G}^{k+1}$ by (14)}.
\FOR {$n=1, \cdots, N$}
\STATE {Update $U^{k+1}_{n}$ by (18)};
\STATE {Update $G^{k+1}_{n}$ by (20)};
\STATE {Update the multiplier $Y^{k+1}_{n}$ by\\
\qquad $Y^{k+1}_{n}=Y^{k}_{n}+\gamma\mu^{k}(\mathcal{G}^{k+1}_{(n)}-G^{k+1}_{n})$.}
\ENDFOR
\STATE {Update $\mu^{k+1}$ by $\mu^{k+1}=\textup{min}(\rho\mu^{k},\,\mu_{max})$.}
\STATE {Check the convergence condition,\\
$\;\;\;\,\textup{max}\{\|\mathcal{G}^{k+1}_{(n)}-G^{k+1}_{n}\|_{F},\;n=1,\cdots, N\}<\textup{Tol}$.}
\ENDWHILE
\OUTPUT $\mathcal{G}$ and $U_{n},\;n=1,\cdots, N$.\\
\end{algorithmic}
\end{algorithm}

\subsection{Complexity Analysis}
We now discuss the time complexity of our NCTD algorithm. For the problem (11), the main running time of our NCTD algorithm is consumed by performing SVD for the proximal operator and some multiplications. The time complexity of performing the proximal operator in (20) is $O_{1}:=O(\sum{{R^{2}_{n}}\Pi_{j\neq{n}}}{R_{j}})$. The time complexity of some multiplication operators is $O_{2}:=O(\sum{{I_{n}R_{n}}{\Pi_{j\neq{n}}}{I_{j}}}+\sum{R_{n}}\Pi_{j\neq{n}}{I_{j}R_{j}})$ and $O_{3}:=O(\sum{R^{2}_{n}I_{n}})$. Thus, the total time complexity of our NCTD method is $O(T(O_{1}+O_{2}+O_{3}))$, where $T$ is the number of iterations.

\section{Experimental Results}
In this section, we evaluate both the effectiveness and efficiency of our methods for solving tensor decomposition problems using both synthetic and real-world data. All experiments were performed on an Intel(R) Core (TM) i5-4570 (3.20 GHz) PC running Windows 7 with 8GB main memory.
\subsection{Synthetic Data}
In this part, we generated a low-rank $N$th-order tensor $\mathcal{T}\in \mathbb{R}^{{I_{1}}\times{I_{2}}\cdots\times{I_{N}}}$, which was used as the ground truth data. The tensor data followed the Tucker model, i.e., $\mathcal{T}=\mathcal{G}\times_{1}U_{1}\cdots\times_{N}U_{N}$, where the core tensor $\mathcal{G}\in\mathbb{R}^{{r}\times{r}\cdots\times{r}}$ and the factor matrices $U_{n}$ were generated with i.i.d. standard Gaussian entries. The order of the tensors varied from three to four, and the rank $r$ was set to $\{5, 10, 20\}$. Finally, we decomposed the input tensor $\mathcal{T}+\delta\Delta$ by our CTD and NCTD methods and the state-of-the-art algorithms including HOSVD (Vannieuwenhoven et al., 2012), HOOI (Lathauwer et al., 2000b), Mixture (Tomioka et al., 2013) and ADMM (Gandy et al., 2011), where $\Delta$ is a noise tensor with independent normally distributed entries, and $\delta$ is set to 0.02.

We set the tensor ranks $R_{n}=\lfloor1.2r\rfloor,\,n=1,\ldots,N$ and $\textup{Tol}=10^{-5}$ for all these algorithms. We set the regularization parameter $\lambda=100$ for our methods. Other parameters of Mixture, ADMM, HOSVD and HOOI are set to their default values. The relative square error (RSE) of the estimated tensor $\mathcal{X}$ is given by $\textup{RSE}=\|\mathcal{X}-\mathcal{T}\|_{F}/\|\mathcal{T}\|_{F}$. The average experimental results (RSE and time cost) of 50 independent runs are shown in Table 1, where the order of tensor data varies from three to four. From the results shown in Table 1, we can see that our methods can yield much more accurate solutions, and outperform their individual competitors, HOSVD, HOOI, Mixture and ADMM, in terms of both RSE and efficiency.
\begin{table*}[t]
\renewcommand{\arraystretch}{1.3}
\scriptsize
\begin{center}
\caption{Performance comparison of estimation accuracy (RSE) and running time (seconds) on the synthetic data:}
(a) Tensor size: $200\times{200}\times{200}$\\
\begin{tabular}[c]{ccccccc|cccccc} 
\hline
\multirow{2}{*}{\ } & \multicolumn{2}{c}{HOSVD} & \multicolumn{2}{c}{HOOI} & \multicolumn{2}{c|}{NCTD} & \multicolumn{2}{c}{ADMM} &  \multicolumn{2}{c}{Mixture} & \multicolumn{2}{c}{CTD}\\
\ Rank      & RSE  & Time 	& RSE 	 & Time	     & RSE 	 & Time                      & RSE 	 & Time   & RSE 	& Time       & RSE         & Time \\
\hline
\ 5  &3.27e-02  &19.46	&3.27e-02	&570.85  	&\textbf{6.52e-03}  &\textbf{16.50}  &1.47e-02  &2946.27	&1.46e-02	&4792.74  &\textbf{6.52e-03}  &\textbf{88.74}\\
\ 10 &3.33e-02  &20.36	&3.32e-02	&468.65	    &\textbf{6.60e-03}  &\textbf{15.25}  &1.48e-02  &2512.85	&1.47e-02	&4567.61  &\textbf{6.59e-03}  &\textbf{85.62}\\
\ 20 &3.34e-02  &21.02	&3.33e-02	&380.87 	&\textbf{6.62e-03}  &\textbf{13.74}  &1.48e-02  &2230.64	&1.47e-02	&4235.65  &\textbf{6.61e-03}  &\textbf{72.69}\\
\hline
\end{tabular}\\
\end{center}
\begin{center}
{(b) Tensor size: $60\times{60}\times{60}\times{60}$\\}
\begin{tabular}[c]{ccccccc|cccccc} 
\hline
\multirow{2}{*}{\ } & \multicolumn{2}{c}{HOSVD} & \multicolumn{2}{c}{HOOI} & \multicolumn{2}{c|}{NCTD} & \multicolumn{2}{c}{ADMM} &  \multicolumn{2}{c}{Mixture} &  \multicolumn{2}{c}{CTD}\\
\ Rank     & RSE   & Time    & RSE  & Time	  & RSE             & Time              & RSE    & Time     & RSE     & Time      & RSE          & Time \\
\hline
\ 5  &3.89e-02	&25.14	 &3.87e-02	&964.26   &\textbf{7.17e-03}  &\,\;\textbf{20.35}\;  &1.64e-02  &3378.06  &1.63e-02	&5175.79   &\textbf{7.16e-03} &\;\textbf{105.01}\;\\
\ 10 &3.91e-02	&23.94	 &3.90e-02	&607.17   &\textbf{6.30e-03}  &\,\;\textbf{18.89}\;  &1.64e-02  &3280.32  &1.64e-02	&4971.15   &\textbf{6.30e-03} &\;\textbf{101.30}\;\\
\ 20 &3.92e-02	&22.81	 &3.91e-02	&415.43   &\textbf{5.67e-03}  &\,\;\textbf{17.63}\;  &1.65e-02  &3031.54  &1.65e-02	&4773.66   &\textbf{5.67e-03} &\;\textbf{98.45}\;\\
\hline
\end{tabular}
\end{center}
\end{table*}

A phase transition plot uses greyscale colors to depict how likely a certain kind of low-rank tensors can be recovered by these algorithms for a range of different given ranks and noise variances $\delta$. Phase transition plots are important means to compare the performance of different tensor estimation methods. If the relative error $\textup{RSE}\leq10^{-2}$, the estimation is regarded as successful. Figure 2 depicts the phase transition plots of HOSVD, HOOI, CTD and NCTD on the third-order tensor data with the rank $r=10$, where the given tensor ranks $R_{n},\,n=1,2,3$ varied from 10 to 50 with increment 4 and $\delta$ from 0 to 0.05 with increment 0.005. For each setting, 50 independent trials were run. From the experimental results shown in Figure 2, we can see that CTD and NCTD perform significantly better than HOSVD and HOOI.
\begin{figure}[!htb]
\centering
\subfigure[HOSVD]{\includegraphics[width=0.45\linewidth]{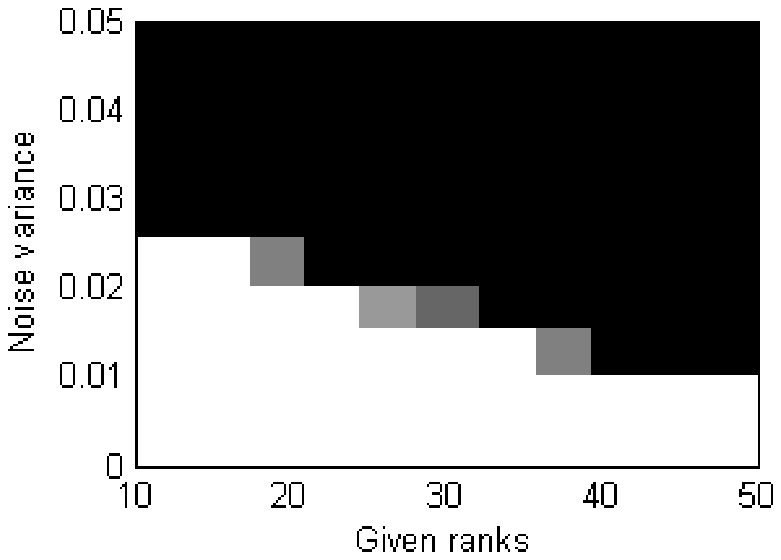}\label{fig_first_case}}
\subfigure[HOOI]{\includegraphics[width=0.45\linewidth]{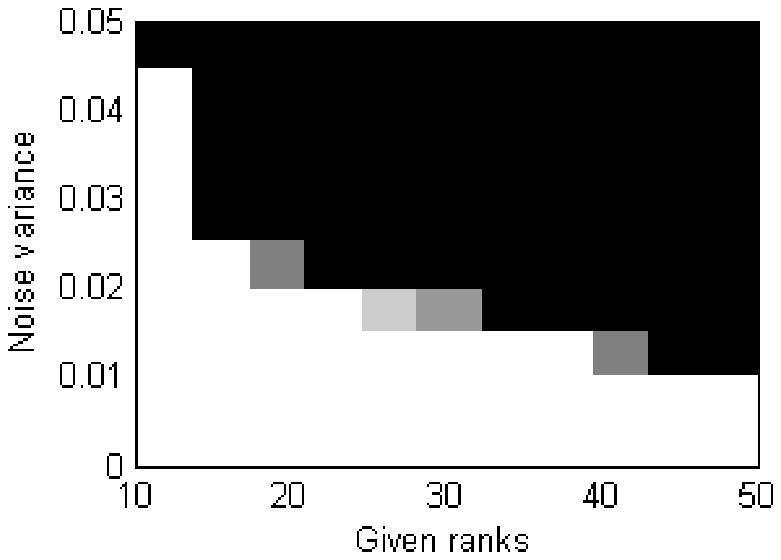}\label{fig_second_case}}
\subfigure[CTD]{\includegraphics[width=0.45\linewidth]{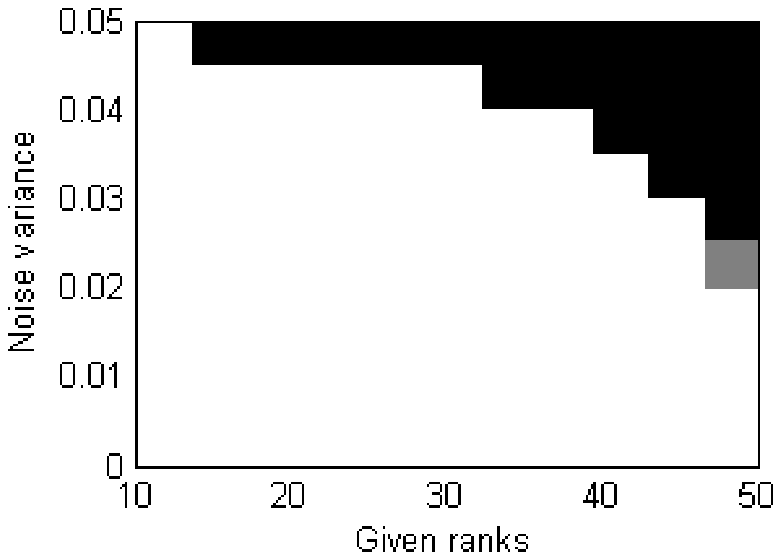}\label{fig_second_case}}
\subfigure[NCTD]{\includegraphics[width=0.45\linewidth]{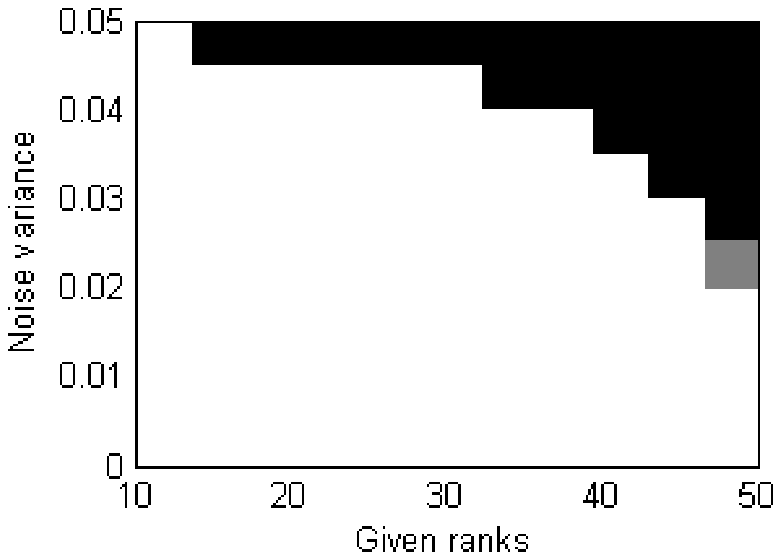}\label{fig_first_case}}
\caption{Phase transition plots for different methods on third-order low-rank tensors with Gaussian random noise, where white denotes perfect estimation in all experiments, and black denotes failure for all experiments.}
\label{fig_sim}
\end{figure}

\subsection{MRI Data}
This part compares our CTD and NCTD methods, HOSVD and HOOI on a $181\times217\times181$ brain MRI data used in (Liu et al., 2009). This data set is approximately low-rank: for the three mode unfoldings, the numbers of singular values larger than 1\% of the largest one are 17, 21, and 17, respectively. Figure 3 shows the average relative errors and running times of ten independent trials for each setting of the given ranks. From the results, we see that our CTD and NCTD methods consistently attain much lower relative errors than those by HOSVD and HOOI. Moreover, our NCTD method is usually faster than the other methods.
\begin{figure}[!htb]
\centering
\subfigure[RSE vs. ranks]{\includegraphics[width=0.46\linewidth]{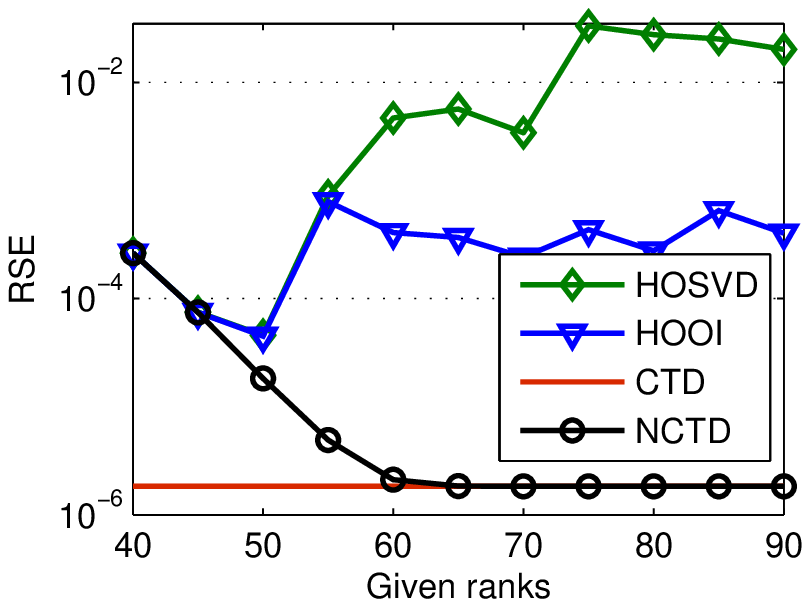}\label{fig_first_case}}
\subfigure[Time vs. ranks]{\includegraphics[width=0.46\linewidth]{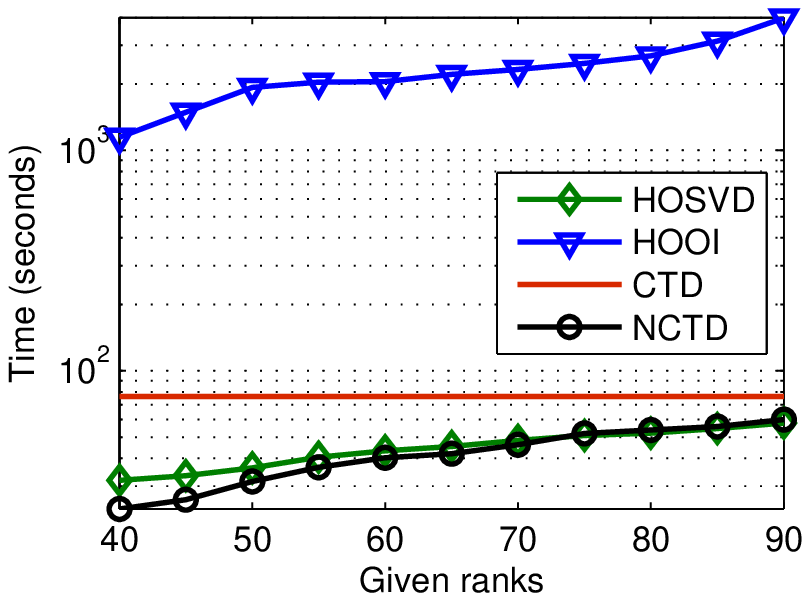}\label{fig_second_case}}
\caption{Comparison of HOSVD, HOOI, CTD and NCTD in terms of estimation accuracy (RSE) and time cost (in the logarithmic scale) on the brain MRI data set.}
\label{fig_sim}
\end{figure}

\section{Conclusions}
In this paper we first proposed a convex trace norm regularized tensor decomposition method, which can automatically determine the number of factors in each mode through our optimization scheme. In addition, by considering the low-rank structure of input tensors, we analyzed the equivalence relationship of the trace norm between a low-rank tensor and its core tensor. Then, we cast the non-convex tensor decomposition model into a weighted combination of multiple much-smaller-scale matrix trace norm minimization. Finally, we developed two efficient parallel ADMM algorithms for solving the proposed problems. Convincing experimental results demonstrate that our regularized formulation is reasonable, and our methods are robust to noise or outliers. Moreover, our tensor decomposition methods can handle some tensor recovery problems, such as tensor completion, and low-rank and sparse tensor decomposition.

\section{Acknowledgements}
We thank the reviewers for their useful comments that have helped improve the paper significantly. The first and third authors are supported by the CUHK Direct Grant No. 4055017 and Microsoft Research Asia Grant 6903555.

\section{References}
\begin{quote}
\begin{small}
Cai, J.; Candes, E.; and Shen, Z. 2010. A singular value thresholding algorithm for matrix completion. \emph{SIAM J. Optim.} 20(4): 1956--1982.

De Lathauwer, L., and Vandewalle, J. 2004. Dimensionality reduction in higher-order signal processing and rank-(r1, r2, . . . , rn) reduction in multilinear algebra. \emph{Linear. Algebra. Appl.} 391: 31--55.

Deng, W.; Lai, M.; and Yin, W. 2013. On the $o(1/k)$ convergence and parallelization of the alternating direction method of multipliers. ArXiv:1312.3040.

Gandy, S.; Recht, B.; and Yamada. I. 2011. Tensor completion and low-n-rank tensor recovery via convex optimization. \emph{Inverse Problem}, 27(2) 025010.

Goldfarb, D., and Qin, Z. 2014. Robust low-rank tensor recovery: models and algorithms. \emph{SIAM J. Matrix Anal. Appl.}

Harshman, R.A. 1970. Foundations of the PARAFAC procedure: models and conditions for an ``explanatory" multi-modal factor analysis. \emph{UCLA working papers in phonetics} 16: 1--84.

He, B. 2009. Parallel splitting augmented Lagrangian methods for monotone structured variational inequalities. \emph{Comput. Optim. Appl.} 42: 195--212.

Hillar, C.J., and Lim, L.H. 2013. Most tensor problems are NP hard. \emph{J. ACM}.

Keshavan, R.H.; Montanari, A.; and Oh, S. 2010. Matrix completion from a few entries. \emph{IEEE Trans. Inform. Theory} 56(6): 2980--2998.

Kolda, T.G., and Bader, B.W. 2009. Tensor decompositions and applications. \emph{SIAM Review} 51(3): 455--500.

Lathauwer, L.; Moor, B.; and Vandewalle, J. 2000a. A multilinear singular value decomposition. \emph{SIAM J. Matrix. Anal. Appl.} 21(4): 1253--1278.

Lathauwer, L.; Moor, B.; and Vandewalle, J. 2000b. On the best rank-1 and rank-($R_{1},R_{2},\cdots,R_{N}$) approximation of high-order tensors. \emph{SIAM J. Matrix. Anal. Appl.} 21(4): 1324--1342.

Lin, Z.; Liu, R.; and Su, Z. 2011. Linearized alternating direction method with adaptive penalty for low-rank representation. In \emph{NIPS}, 612--620.

Liu, J.; Musialski, P.; Wonka, P.; and Ye, J. 2009. Tensor completion for estimating missing values in visual data. In \emph{ICCV} 2114--2121.

Liu, Y.; Shang, F.; Cheng, H.; Cheng, J.; and Tong, H. 2014. Factor matrix trace norm minimization for low-rank tensor completion. In \emph{SDM}.

Morup, M. 2011. Applications of tensors (multiway array) factorizations and decompositions in data mining. \emph{WIREs: Data Min. Knowl. Disc.} 1(1): 24--40.

Narita, A.; Hayashi, K.; Tomioka, R.; and Kashima, H. 2012. Tensor factorization using auxiliary information. \emph{Data Min. Knowl. Disc.} 25(2): 298--324.

Nick, H. 1995. Matrix procrustes problems.

Shi, Z.; Han, J.; Zheng, T.; and Li, J. 2013. Guarantees of augmented trace norm models in tensor recovery. In \emph{IJCAI} 1670--1676.

Shang, F.; Jiao, L.C.; Shi, J.; and Chai, J. 2011. Robust positive semidefinite L-Isomap ensemble. \emph{Pattern Recogn. Lett.} 32(4): 640--649.

Signoretto, M.; Dinh, Q.T.; Lathauwer, L.D.; and Suykens, J. 2014. Learning with tensors: a framework based on convex optimization and spectral regulation. \emph{Mach. Learn.}

Tomioka, R.; Suzuki, T.; Hayashi, K.; and Kashima, H. 2011. Statistical Performance of Convex Tensor Decomposition. In \emph{NIPS} 972--980.

Tomioka, R., and Suzuki, T. 2013. Convex tensor decomposition via structured Schatten norm regularization. In \emph{NIPS} 1331--1339.

Tucker, L.R. 1966. Some mathematical notes on three-mode factor analysis. \emph{Psychometrika} 31: 279--311.

Vannieuwenhoven, N.; Vandebril, R.; and Meerbergen, K. 2012. A new truncation strategy for the higher-order singular value decomposition. \emph{SIAM J. Sci. Comput.} 34(2): 1027--1052.

Wen, Z.; Yin, W.; and Zhang, Y. 2012. Solving a low-rank factorization model for matrix completion by a nonlinear successive over-relaxation algorithm. \emph{Math. Prog. Comp.} 4(4): 333--361.

Yilmaz, Y.K.; Cemgil, A.T.; and Simsekli, U. 2011. Generalized coupled tensor factorization. In \emph{NIPS} 2151--2159.
\end{small}
\end{quote}

\balance

\section*{Convergence Behaviors of Our Algorithms}
We also study the convergence behaviors of our CTD and NCTD algorithms on the synthetic tensor data of size $200\times200\times200$ with the given ranks, $R_{1}=R_{2}=R_{3}=20$, as shown in Figure 1, where the ordinate is the residual of $\max\{\|\mathcal{G}^{k}_{(n)}-G^{k}_{n}\|_{F}, n=1,\cdots,N\}$ or $\|\mathcal{X}^{k}-\mathcal{T}\|_{F}$, or the relative change of $\mathcal{X}^{k}$, and the abscissa denotes the number of iterations. Moreover, we also provide the convergence results of HOOI. We can observe that the relative change or the residual of CTD and NCTD algorithms drops much quickly, and converges much fast within 50 iterations. Especially, the relative change or the residual of CTD and NCTD drops much more quickly than HOOI.

\begin{figure}[!ht]
\centering
\includegraphics[width=0.494\linewidth]{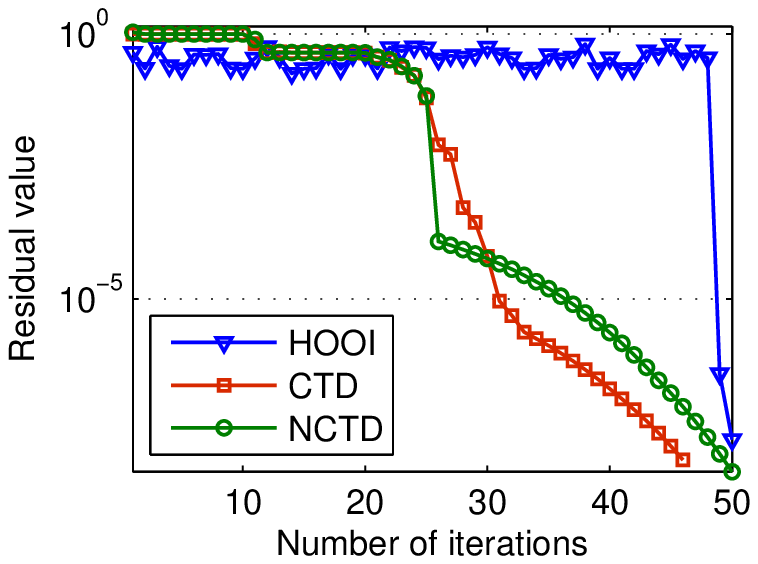}\label{fig_first_case}
\includegraphics[width=0.494\linewidth]{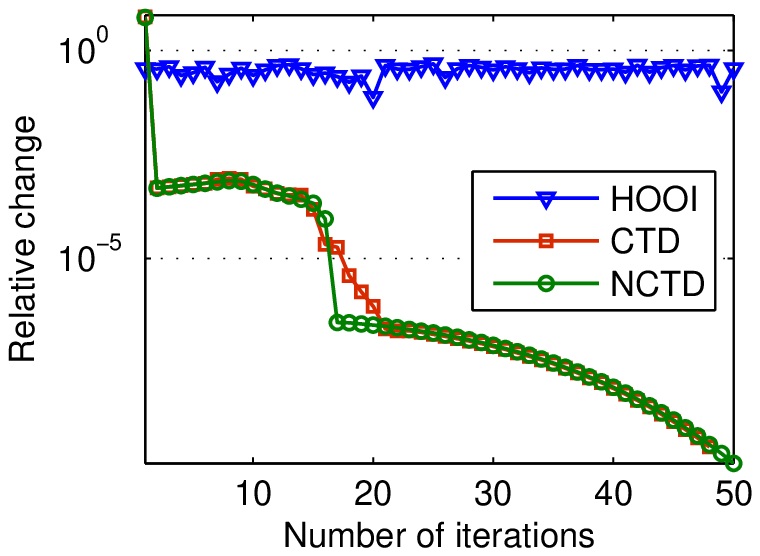}\label{fig_second_case}
\centering
\caption{Convergence behaviors of HOOI, our CTD and NCTD algorithms on the synthetic tensor data of size $200\times200\times200$. Left: the residual of $\max\{\|\mathcal{G}^{k}_{(n)}-G^{k}_{n}\|_{F}\}$ or $\|\mathcal{X}^{k}-\mathcal{T}\|_{F}$. Right: the relative change of $\mathcal{X}^{k}$.}
\label{fig_sim}
\end{figure}

\subsection*{Rank Estimation}
In this part, we evaluate the ability of our CTD method to estimate the tensor ranks, as shown in Figure 2. As in the experimental section, we randomly generated $\mathcal{T}$ of size $200\times200\times200$ with the tensor ranks $r_{1}=r_{2}=r_{3}=r=$10 or 20. From Figure 2, we see that our CTD method can accurately estimate the rank of each mode unfolding of the tensor.

\begin{figure}[!ht]
\centering
\subfigure[True rank $r = 10$]{\includegraphics[width=0.494\linewidth]{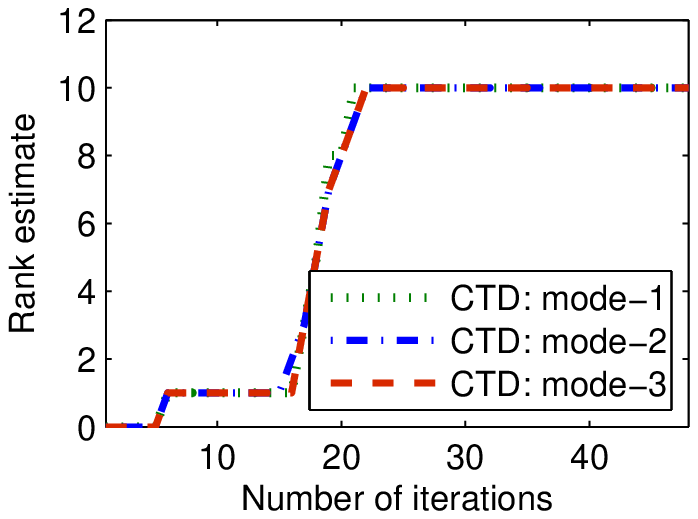}\label{fig_first_case}}
\subfigure[True rank $r = 20$]{\includegraphics[width=0.494\linewidth]{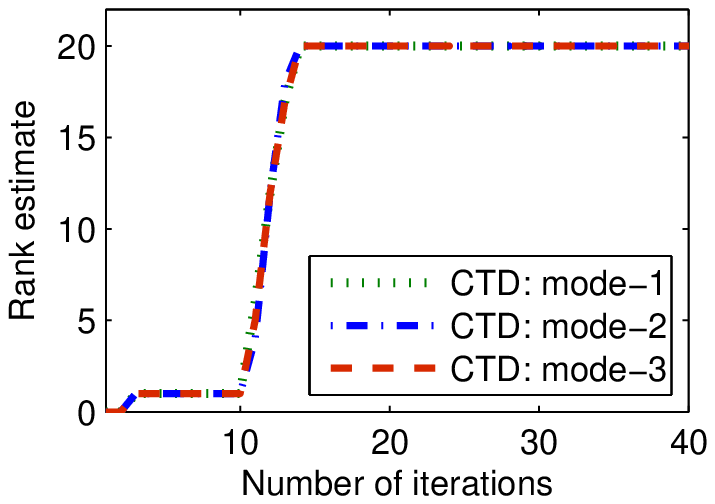}\label{fig_second_case}}
\centering
\caption{Estimation of the rank of each mode unfolding on third-order low-rank tensors.}
\label{fig_sim}
\end{figure}

\section*{Robustness Against Outliers}
Figure 3 depicts the phase transition plots of HOSVD, HOOI, CTD and NCTD on the third-order tensor data with rank $r=10$, where the given tensor ranks $R_{n},\,n=1,2,3$ varied from 10 to 50 with increment 4 and the error sparsity ratio $\|\Delta\|_{0}/\Pi_{n}I_{n}$ from 0 to 0.05 with increment 0.005. We generated $\Delta$ as a sparse tensor whose non-zero entries are independent and uniformly distributed in the range $[-1,1]$, and whose support is chosen uniformly at random. For each setting, 50 independent trials were run. From the experimental results shown in Figure 3, it is clear that our methods, CTD and NCTD, significantly outperform HOSVD and HOOI.

\begin{figure}[!ht]
\centering
\subfigure[HOSVD]{\includegraphics[width=0.47\linewidth]{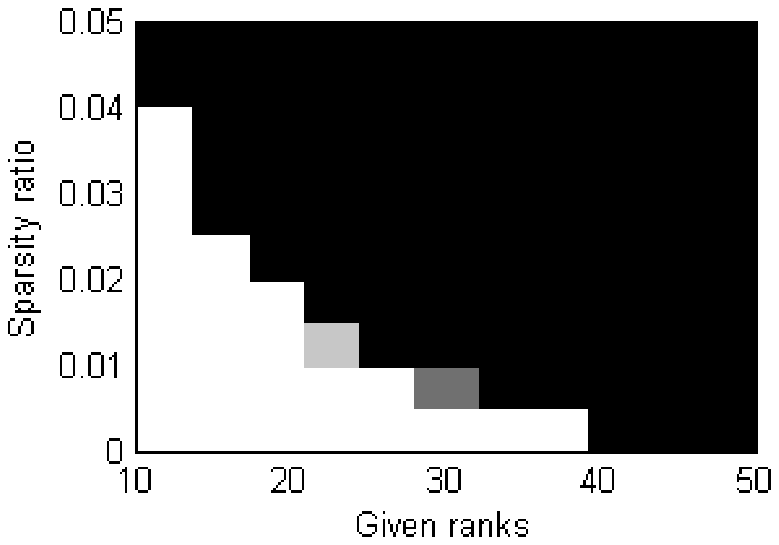}\label{fig_first_case}}
\subfigure[HOOI]{\includegraphics[width=0.47\linewidth]{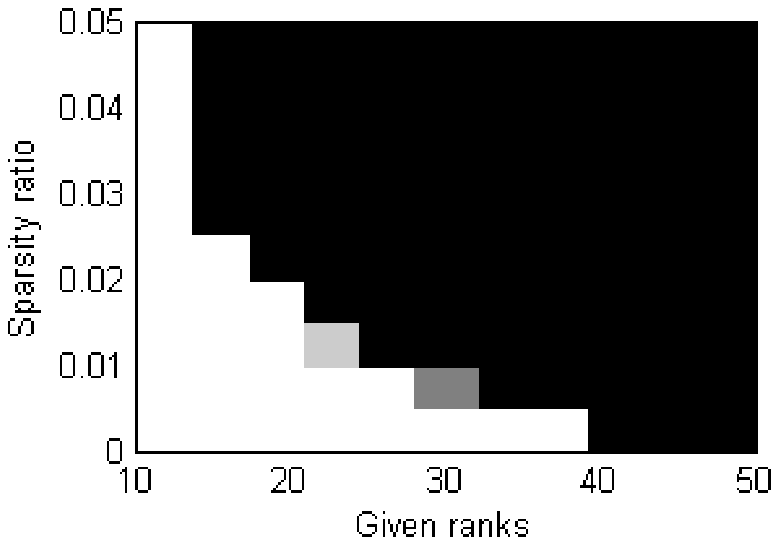}\label{fig_second_case}}
\subfigure[CTD]{\includegraphics[width=0.47\linewidth]{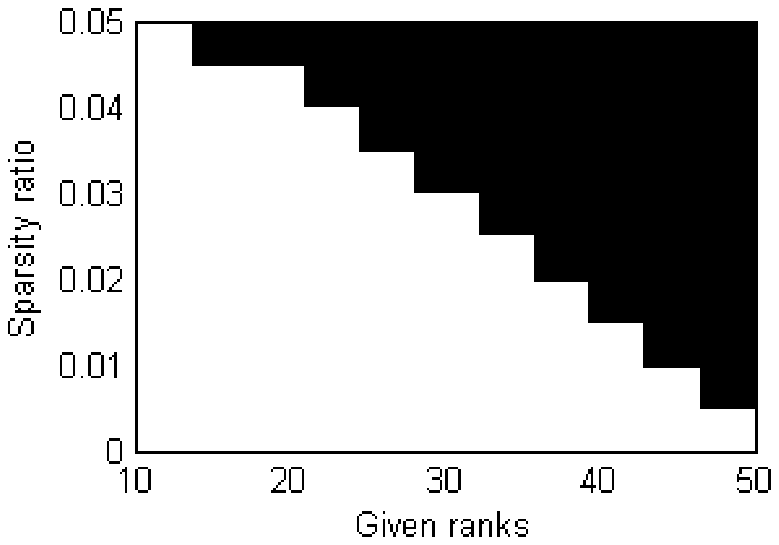}\label{fig_second_case}}
\subfigure[NCTD]{\includegraphics[width=0.47\linewidth]{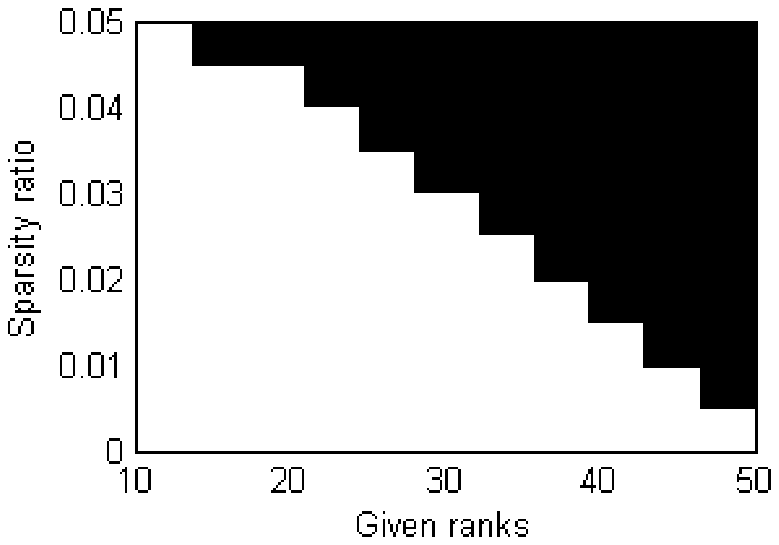}\label{fig_first_case}}
\caption{Phase transition plots for different methods on third-order low-rank tensors with outliers, where white denotes perfect estimation in all experiments, and black denotes failure for all experiments.}
\label{fig_sim}
\end{figure}

\section*{APPENDIX A}
\textbf{\emph{Proof of Theorem 2}}:
Before giving a proof of Theorem 2, we will first present some properties of matrices and tensors in the following.

\textbf{Property 1}: Let $U\in \mathbb{R}^{m\times p}$, $S\in \mathbb{R}^{p\times q}$, and $V\in \mathbb{R}^{n\times q}$, then
\begin{displaymath}
\|USV^{T}\|_{\textup{tr}}=\|S\|_{\textup{tr}},
\end{displaymath}
where $U^{T}U=I_{p}$ and $V^{T}V=I_{q}$.

\textbf{Property 2}: Let $A\in \mathbb{R}^{m\times n}$, $B\in \mathbb{R}^{p\times q}$, and $C$ and $D$ are two matrices of suitable sizes, then we have the following results:\\
(I). $(A\otimes B)\otimes C=A\otimes (B\otimes C)$.\\
(II). $(A\otimes B)(C\otimes D)=AC\otimes BD$.\\
(III). $(A\otimes B)^{T}=A^{T}\otimes B^{T}$.

\textbf{Property 3}: Let $\mathcal{X}=\mathcal{G}\times_{1}U_{1}\times_{2}\cdots\times_{N}U_{N}$, where $\mathcal{X}\in \mathbb{R}^{{I_{1}}\times{I_{2}}\cdots\times{I_{N}}}$, $\mathcal{G}\in \mathbb{R}^{{R_{1}}\times{R_{2}}\cdots\times{R_{N}}}$, and $U_{n}\in \mathbb{R}^{I_{n}\times R_{n}},n=1,\ldots,N$, then
\begin{displaymath}
\mathcal{X}_{(n)}=U_{n}\mathcal{G}_{(n)}(U_{N}\otimes\ldots U_{n+1}\otimes U_{n-1}\ldots \otimes U_{1})^{T}.
\end{displaymath}
\\
\emph{\textbf{Proof of Theorem 2}}: According to Property 3, we have
\begin{displaymath}
\|\mathcal{X}_{(n)}\|_{\textup{tr}}=\|U_{n}\mathcal{G}_{(n)}(U_{N}\otimes\ldots U_{n+1}\otimes U_{n-1}\ldots \otimes U_{1})^{T}\|_{\textup{tr}}.
\end{displaymath}

Let $\emph{P}_{n}=U_{N}\otimes\ldots U_{n+1}\otimes U_{n-1}\ldots \otimes U_{1}$ and $U^{T}_{n}U_{n}=I_{R_{n}},n=1,\ldots,N$, and according to Properties 2 and 3, we have the following conclusion:
\begin{displaymath}
\begin{split}
&P^{T}_{n}P_{n}\\
=&(U_{N}\otimes\ldots \otimes U_{n+1}\otimes U_{n-1}\otimes\ldots \otimes U_{1})^{T}(U_{N}\otimes\ldots\otimes\\
& U_{n+1}\otimes U_{n-1}\otimes\ldots\otimes U_{1}),\\
=&(U^{T}_{N}\otimes\ldots\otimes U^{T}_{n+1}\otimes U^{T}_{n-1}\otimes\ldots \otimes U^{T}_{1})(U_{N}\otimes\ldots\otimes\\
&U_{n+1}\otimes U_{n-1}\otimes\ldots \otimes U_{1}),\\
=&(U^{T}_{N}\otimes\ldots\otimes U^{T}_{n+1}\otimes U^{T}_{n-1}\otimes\ldots \otimes U^{T}_{2})(U_{N}\otimes\\
&\ldots\otimes U_{n+1}\otimes U_{n-1}\otimes\ldots \otimes U_{2})\otimes(U^{T}_{1}U_{1}),\\
=&(U^{T}_{N}\otimes\ldots\otimes U^{T}_{n+1}\otimes U^{T}_{n-1}\otimes\ldots \otimes U^{T}_{2})(U_{N}\otimes\\
&\ldots\otimes U_{n+1}\otimes U_{n-1}\otimes\ldots \otimes U_{2})\otimes I_{1},\\
\vdots\\
=&I_{N}\otimes \ldots I_{n+1}\otimes I_{n-1} \otimes \ldots I_{2}\otimes I_{1},\\
=&\widetilde{I}_{n},
\end{split}
\end{displaymath}
where $I_{n}\in \mathbb{R}^{R_{n}\times R_{n}},n=1,\ldots,N,\,\widetilde{I}_{n}\in\mathbb{R}^{J_{n}\times J_{n}}$, and $ J_{n}=\Pi_{j\neq n}R_{j}$.

According to Property 1, and $P^{T}_{n}P_{n}=\widetilde{I}_{n}$, we have
\begin{displaymath}
\begin{split}
\|\mathcal{X}_{(n)}\|_{\textup{tr}}&=\|U_{n}\mathcal{G}_{(n)}(U_{N}\otimes\ldots U_{n+1}\otimes U_{n-1}\ldots \otimes U_{1})^{T}\|_{\textup{tr}}\\
&=\|\mathcal{G}_{(n)}\|_{\textup{tr}}.
\end{split}
\end{displaymath}
This completes the proof. \;\;\quad\quad\quad\quad\quad\quad\quad\quad\quad\quad\quad\quad$\Box$

\section*{APPENDIX B}
\textbf{\emph{Proof of Theorem 3}}:\\
The optimization problem (13) with respect to $\mathcal{G}$ is written by
\begin{equation*}\leqno{(21)\;\;}
\begin{split}
\min_{\mathcal{G}}\,J(\mathcal{G})=\sum^{N}_{n=1}\frac{\mu^{k}}{2}\|\mathcal{G}_{(n)}-G^{k}_{n}+Y^{k}_{n}/\mu^{k}\|^{2}_{F}\\
+\frac{\lambda}{2}\|\mathcal{T}-\mathcal{G}\times_{1}U^{k}_{1}\cdots\times_{N}U^{k}_{N}\|^{2}_{F}.
\end{split}
\end{equation*}
The problem (21) above is a smooth convex optimization problem, thus we can obtain the derivative of the function $J$ in the following form:
\begin{equation*}\leqno{(22)}
\begin{split}
\frac{\partial J}{\partial\mathcal{G}}=&\lambda(\mathcal{G}-\mathcal{T}\times_{1}(U^{k}_{1})^{T}\cdots\times_{N}(U^{k}_{N})^{T})\\
&+\sum^{N}_{n=1}\mu^{k}(\mathcal{G}-\textup{refold}(G^{k}_{n}-Y^{k}_{n}/\mu^{k}))\\
=&(N\mu^{k}+\lambda)\mathcal{G}-\mu^{k}\sum^{N}_{n=1}\textup{refold}(G^{k}_{n}-Y^{k}_{n}/\mu^{k})\\
&-\lambda\mathcal{T}\times_{1}(U^{k}_{1})^{T}\cdots\times_{N}(U^{k}_{N})^{T},
\end{split}
\end{equation*}
where $\textup{refold}(\cdot)$ denotes the refolding of the matrix into a tensor.

Let $\frac{\partial J}{\partial\mathcal{G}}$=0, and
\begin{equation*}\leqno{(23)\quad\quad}
\begin{split}
\mathcal{M}&=\mathcal{T}\times_{1}(U^{k}_{1})^{T}\cdots\times_{N}(U^{k}_{N})^{T},\\ \mathcal{N}&=\sum^{N}_{n=1}\textup{refold}(G^{k}_{n}-Y^{k}_{n}/\mu^{k}),
\end{split}
\end{equation*}
the optimal solution to (21) is given by
\begin{equation*}\leqno{(24)}
\begin{split}
\mathcal{G}^{\ast}=&\frac{\lambda}{\lambda+N\mu^{k}}\mathcal{T}\times_{1}(U^{k}_{1})^{T}\cdots\times_{N}(U^{k}_{N})^{T}\\
&+\frac{\mu^{k}}{\lambda+N\mu^{k}}\sum^{N}_{n=1}\textup{refold}(G^{k}_{n}-Y^{k}_{n}/\mu^{k}),\\
=&\frac{\lambda}{\lambda+N\mu^{k}}\mathcal{M}+\frac{\mu^{k}}{\lambda+N\mu^{k}}\mathcal{N}.
\end{split}
\end{equation*}
This completes the proof. \;\;\quad\quad\quad\quad\quad\quad\quad\quad\quad\quad\quad\quad$\Box$

\section*{APPENDIX C}
\textbf{\emph{Proof of Theorem 4}}:\\
Let
\begin{equation*}\leqno{(25)\;\;}
\begin{split}
H(\mathcal{G},\{U_{n}\})=\sum^{N}_{n=1}\mu^{k}\|\mathcal{G}_{(n)}-G^{k}_{n}+Y^{k}_{n}/\mu^{k}\|^{2}_{F}\\
+\lambda\|\mathcal{T}-\mathcal{G}\times_{1}U_{1}\cdots\times_{N}U_{N}\|^{2}_{F},\\
\end{split}
\end{equation*}
then the closed-form solution of (25) with respect to $\mathcal{G}$ can be obtained by (14), and it can be rewritten as
\begin{equation*}
\mathcal{G}^{\ast}=\frac{1}{\lambda+N\mu^{k}}(\lambda\mathcal{M}+\mu^{k}\mathcal{N}).
\end{equation*}

Hence, we have
\begin{equation*}\leqno{(26)\;\;}
\begin{split}
&H(\mathcal{G}^{\ast},\{U_{n}\})\\
=&\sum^{N}_{n=1}\mu^{k}\|\mathcal{G}^{\ast}_{(n)}-G^{k}_{n}+Y^{k}_{n}/\mu^{k}\|^{2}_{F}+\lambda\|\mathcal{T}\|^{2}_{F}+\lambda\|\mathcal{G}^{\ast}\|^{2}_{F}\\
&-2\lambda\langle\mathcal{T},\mathcal{G}^{\ast}\times_{1}U_{1}\cdots\times_{N}U_{N}\rangle,\\
=&\eta-2\lambda h(U_{1},U_{2},\ldots,U_{N}),
\end{split}
\end{equation*}
where $\eta=\sum^{N}_{n=1}\mu^{k}\|\mathcal{G}^{\ast}_{(n)}-G^{k}_{n}+Y^{k}_{n}/\mu^{k}\|^{2}_{F}+\lambda\|\mathcal{T}\|^{2}_{F}+\lambda\|\mathcal{G}^{\ast}\|^{2}_{F}$ is a constant, and
\begin{displaymath}
h(U_{1},U_{2},\ldots,U_{N})=\langle\mathcal{T},\mathcal{G}^{\ast}\times_{1}U_{1}\cdots\times_{N}U_{N}\rangle.
\end{displaymath}
Combination with the results above proves the theorem.\\
This completes the proof. \;\;\quad\quad\quad\quad\quad\quad\quad\quad\quad\quad\quad\quad$\Box$

\end{document}